\documentclass[11pt]{article}
\usepackage{amsmath}
\usepackage{graphicx}
\usepackage{amsfonts}
\usepackage{amssymb}
\newcommand{\ds}{\displaystyle}

\begin{document}

\title{{Color Visualization of Blaschke Self-Mappings 
of the Real
Projective Plan}}

\author{Cristina Ballantine and
Dorin Ghisa}

\maketitle

\begin{abstract}The real projective plan $P^{2}$ can be endowed with a
dianalytic structure making it into a non orientable Klein surface.
Dianalytic self-mappings of that surface are \textit{projections }of
analytic self-mappings of the Riemann sphere\textbf{\ }$\widehat{\mathbb{C}}$. It is
known that the only analytic bijective self-mappings of $\widehat{\mathbb{C}}$ are
the M\"{o}bius transformations. The Blaschke products are obtained by
multiplying particular M\"{o}bius transformations. They are no longer
one-to-one mappings. However, some of these products can be  \textit{projected%
} on $P^{2}$ and they become dianalytic self-mappings of $P^{2}.$
More exactly, they represent canonical projections of non orientable
branched covering Klein surfaces over $P^{2}.$ This article is devoted to 
color visualization of such mappings. The working tool is the technique of 
\textit{simultaneous continuation} we introduced in previous papers.
Additional graphics and animations are provided on the web site of the project \cite{cri}.

\end{abstract}


\noindent AMS Subject Classification,\emph{\ Primary: 30D50, Secondary: 30F50}

\section{Blaschke Products and Blaschke
Quotients}

The building blocks of Blaschke products are M\"{o}bius transformations of
the form \begin{equation}b_{k}(z)=e^{i\theta _{k}}\frac{a_{k}-z}{%
1-\overline{a}_{k}z},\end{equation} where  $a_{k}\in D:=\{z\in \mathbb{C}\ |\ |z|<1\},$ and $\theta _{k}\in \mathbb{R}.$ We call
them \textit{Blaschke factors.} A finite (infinite) Blaschke product has the
form\begin{equation}w=B(z)=\prod\limits_{k=1}^{n}b_{k}(z),\end{equation}
where \ $n\in \mathbb{N}$, (respectively $n=\infty ).$ In the infinite case it is
customary to take $\ds e^{i\theta _{k}}=\overline{a}_{k}/|a_{k}|.$ We will
adopt these values throughout in this paper. Finite Blaschke products are
meromorphic functions in $\widehat{\mathbb{C}},$ having all the poles outside the
unit disk. Infinite Blaschke products cannot be defined on the set $E$ of
accumulation points of the \textit{zeros} $a_{k}$ of $B,$ yet they are
meromorphic functions in $\widehat{\mathbb{C}}\backslash E.$ If we allow values
of $a_{k}$ with $|a_{k}|$ $>1$ in $(1),$ then $(2)$ becomes the quotient of
two legitimate Blaschke products (see \cite{4}), called a \textit{Blaschke
quotient. }

It is known that the necessary and sufficient condition for an infinite
Blaschke product $(2)$ to be convergent is that $\ds \sum\limits_{n=1}^{\infty
}(1-|a_{k}|)<\infty .$ Moreover, if this condition holds, then the
product converges uniformly on every compact subset of $\mathbb{C}\backslash 
(A\cup E),$ where $A=\{1/\overline{a}_{n}:n\in \mathbb{N}\}$ is the set of \textit{%
poles} of $B.$ It is an easy exercise to show that the same condition is
necessary and sufficient for the convergence of a Blaschke quotient $B.$ The
set $A$ is\ still the set of poles of $B$, some of which might be now
situated inside the unit disc. Again, if the condition holds, then the
convergence is uniform on every compact subset of $\mathbb{C} \backslash (A\cup
E),$ where $E\subset \partial D$ is the set of accumulation points of the
sequence $(a_{k}).$

The papers \cite{1}-\cite{3} deal with the study of global geometric properties of
Blaschke products, while \cite{4} and \cite{5} refer to Blaschke self mappings of the
real projective plan. In the case of infinite Blaschke products and
quotients we restricted our study to the situation where $E$ is a
(generalized) Cantor subset of the unit circle. Such a set is the union of a
discrete non-empty set and a ternary Cantor subset on the unit circle, which
might be the empty set. Obviously, it contains no arc of the unit circle.
This restriction is neither necessary nor sufficient for the convergence of
infinite Blaschke products or quotients, yet it is inclusive enough to allow
the treatment of a wide class of such functions. Most of the results we
obtained in these papers are expected to be true in some different settings.

\section{Blaschke Products and Quotients
Commuting with $h$}

The function $h:\widehat{\mathbb{C}}\rightarrow\widehat{\mathbb{C}}$ defined by $h(z)=-1/\overline{z},$ if $%
z\notin \{0,\infty \},$\ and $h(0)=\infty ,$ $h(\infty )=0,$ is a fixed
point free antianalytic involution. This involution and the functions
commuting with it play a crucial role in the theory of Klein surfaces. It is
therefore interesting to study Blaschke products and quotients commuting
with $h.$ If $B$ is a Blaschke product or quotient, then $B\circ h=h\circ B$
if and only if $B$ is of the form (see \cite{5}):\begin{equation}B(z)=z^{2p+1}\prod_{k=1}^{n\leq \infty }\frac{\overline{a}_{k}%
}{a_{k}}\frac{a_{k}^{2}-z^{2}}{1-\overline{a}_{k}^{2}z^{2}},\end{equation}
where $p$ is an integer and $|a_{k}|$ $\neq 1$.

We are looking for classes of Blaschke products of the types studied in \cite{1}
which also have the property of commuting with $h.$ If  in
formula $(7)$ from \cite{1} we take $a_{1}=-a_{2}=a$, let $n$  be an odd positive
integer and multiply by the factor $z^{n},$ we obtain such a product. It
is of the form:
\begin{equation}B(z)=B_{a}(z)=z^{n}\left[\frac{\overline{a}}{a}\frac{a^{2}-z^{2}}{%
1-\overline{a}^{2}z^{2}}\right]^{n}.\end{equation}

If we let $a$ in $(4)$ take a value with $|a|$ $>1,$ then after the
substitution $a=1/\overline{c},$ we obtain $B_{a}(z)=z^{n}/
B_{c}(z),$ where $|c|<1$. Therefore, $B_{a}$ appears as the quotient of
two legitimate Blaschke products. A similar result is obtained if we allow
in $(3)$ negative values for $p.$ When dealing with Blaschke self-mappings
of the real projective plan, the fact that $|a_{k}|$ in $(3)$ are greater or
less than $1,$ or $p$ is positive or negative does not matter, since $a_{k}$
and $ -1/\overline{a}_{k}$ , as well as $z$ and $ -1/\overline{z}$ are
identified in order to obtain corresponding points on $P^{2}.$ In other
words, in such an instance Blaschke products and quotients must be treated
as a unique class of functions.
\bigskip

\textbf{Theorem 2.1.} \textit{The fundamental domains of} $(\widehat{\mathbb{C}},B),$
 \textit{where} $B$ \textit{is given by} $(4)$ \textit{are bounded by
consecutive arcs connecting} $a$ \textit{and} $1/\overline{a},$
$-a$ \textit{and} $-1/\overline{a}$, $b$ \textit{and} $1/\overline{b},$ \textit{respectively} $-b$ \textit{and} $-1/\overline{b},$
\textit{a part of the line determined by} $0$ \textit{and} $a,$ \textit{as well as its symmetric with respect to the unit circle, and some infinite
rays issued from the origin. Here} $b$ \textit{is any one of the four
solutions of the equation} $B^{\prime }(z)=0$ \textit{which is not a zero of} $B.$ \textit{Every fundamental domain is mapped conformally by} $B$ \textit{on the} $w$\textit{-plane from which a part of the real axis has been removed.}\bigskip

\textbf{Proof:\ }If $a=re^{i\alpha },$ then the equation $B(z)=\tau
e^{in\alpha },$ $\tau \geq 0$ is equivalent to\bigskip\ the set of equations
\begin{equation}z\frac{r^{2}-(e^{-i\alpha }z)^{2}}{1-(e^{-i\alpha }rz)^{2}}=\tau
^{1/n}e^{i\alpha }\omega _{k},\qquad \tau \geq 0,\qquad k=0,1,...,n-1,\end{equation}
where $\omega _{k}$ are the n-th roots of unity. These equations become,
after the substitution $u=e^{-i\alpha }z$,\begin{equation}u\frac{r^{2}-u^{2}}{1-r^{2}u^{2}}=\tau ^{1/n}\omega
_{k},\qquad \tau \geq 0,\qquad k=0,1,...,n-1.\end{equation}
Each one of the equations $(5)$ is a third degree equation and therefore has
three roots (distinct or not). It is known (see \cite{2}) that for $|\zeta |=1$ we have $B^{\prime }(\zeta )\neq 0,$ and also that $|B(z)|=1$ if and
only if $|z|=1.$ Therefore, for $\tau =1,$ the solutions of the
equations $(5)$ are distinct and they are all points on the unit circle. We
notice that $\zeta _{0}=-e^{i\alpha }$ is a solution of the equation $(5)$
corresponding to $k=0$ and that the other solutions of $(5)$ are two by two
symmetric with respect to the line $z=te^{i\alpha },$ $t \in \mathbb{R}.$ Let us
denote them by $\zeta _{1},\zeta _{2},...,\zeta _{3n-1}$ visited
counter-clockwise on the unit circle. The simultaneous continuation in the
unit disc over $w(\tau )=(1-\tau) e^{in\alpha },$ $0\leq \tau \leq 1$, starting
from these points produces $3n$ arcs $z_{j}(\tau )$ converging $n$ of them
to $a,$ $n$ of them to $-a$ and $n$ of them to $0$ and arriving there when $\tau =1.$ Indeed, the points $\pm a,$ and $0$ as multiple zeros of $B$ of
order of multiplicity $n,$ are branch points of order $n$ of the branched
covering surface $(\widehat{\mathbb{C}},B)$ (see \cite{2} and \cite{3}). Due to the fact
that $B$ and $h$ commute, so are $\pm 1/\overline{a}$ and $\infty .$ These
last points can be reached by simultaneous continuation over $w(\tau )=\tau
e^{in\alpha },$ $\tau \in (1,+\infty ).$ The other branch points can be
found by solving the equation $\ds\frac{d}{du}\left[u\frac{r^{2}-u^{2}}{1-r^{2}u^{2}}\right]=0,$ where $u=e^{-i\alpha }z.$ The solutions of this equation are, in terms of $z:$
\begin{equation}\pm \frac{1}{r\sqrt{2}}\sqrt{3-r^{4}\pm \sqrt{(3-r^{4})^{2}-4r^{4}}}e^{i\alpha }.\end{equation}
It is obvious that $(3-r^{4})^{2}-4r^{4}>0$ and $3-r^{4}-\sqrt{(3-r^{4})^{2}-4r^{4}}>0$. Hence, these branch points are all on the line $%
z(t)=e^{i\alpha }t,$ $t\in \mathbb{R}$ passing through $\pm a$ and $0,$ as expected
(see \cite{7}). In fact, denoting by $b$ the solution situated between $0$ and $%
a,$ i.e. \begin{equation}b=\frac{1}{r\sqrt{2}}\sqrt{3-r^{4}-\sqrt{(3-r^{4})^{2}-4r^{4}%
}}e^{i\alpha},\end{equation}
the other solutions are $-b,$ $1/\overline{b},$ and $-1/\overline{b}.$
We notice that
$$B(e^{i\alpha }t)=e^{in\alpha }t^{n}[(t^{2}-r^{2})/
(1-r^{2}t^{2})]^{n}$$
and since $|b|<r$, as it can be easily checked, this formula shows that
the branch points $\pm b$ are reached by $z(t),$ $t<1$, after passing through 
$\pm a,$ which correspond to $\tau =1.$ This means that we need to extend
the simultaneous continuation for negative values of $1-\tau ,$ more exactly
for values of $1-\tau $ between $0$ and $-|B(b)|$ to let these arcs meet each
other. Similarly, the points $\pm 1/\overline{b}$ are reached by $z(t),$ $t >1$, after passing through $\pm 1/\overline{a},$ which correspond to $\tau
=+\infty $. Again, we need to let $\tau $ vary through negative values
less than $-|B(1/\overline{b})|.$ The consecutive arcs obtained above border domains which are mapped by $B$ conformally on the $w$%
-plane from which the ray $w(t)=te^{in\alpha },t>0$, has been removed
(fundamental domains). The mappings are continuous also on the border of
every fundamental domain and there is a continuous passage from the mapping
of one fundamental domain to that of any adjacent one. Moreover, $B$ is
conformal also on the boundaries of the fundamental domains, except for the
branch points. We notice that the segment between $-a$ and $-1/\overline{a}$
is on the border of two fundamental domains of $B$, since $-e^{i\alpha }$ 
is always a solution of the equation $(4)$ corresponding to $k=0,$ while the
open interval between $a$ and $1/\overline{a}$ is always inside a
fundamental domain of $B$. There are $n$ fundamental domains bounded by arcs
connecting $a$ and $1/\overline{a},$ as well as $n$ fundamental domains
bounded by arcs connecting $-a$ and $-1/\overline{a},$ and $n$ unbounded
fundamental domains. Two of these last domains contain on their boundary the
segment between $0\ $\ and $a$ and the segment between$1/\overline{a}$ and $%
\infty $ on the ray $z(t)=e^{i\alpha }t,$ $t>0,$ and other two the segment
between $0$ and $-a$ and the segment between $-1/\overline{a}$ and $\infty $
on the ray $z(t)=e^{i\alpha }t,$ $t<0.$ The other unbounded fundamental
domains have as boundaries couples of arcs extending from $0$ to $\infty .$
\bigskip

\textbf{Theorem 2.2: }\textit{Let us denote by }$\Omega _{j}$ \textit{the
fundamental domain containing the arc of the unit circle between }$\zeta
_{j} $\textit{\ and }$\zeta _{j+1},$\textit{\ }$j=0,1,...,3n-1$\textit{\ and
let}
$$S_{k}(z)=B_{\Omega _{(k+j)\!\!\! \!\! \pmod{3n}}}^{-1}\circ B_{\Omega
_{j}}(z),\qquad j=0,1,...,3n-1,\qquad k=0,1,...,3n-1.$$
\textit{Then} $\{S_{k}\}$ \textit{is a group under the composition law }$%
S_{k}\circ S_{j}=S_{(k+j)\!\! \pmod{3n}},$\textit{\ where }$S_{0}$\textit{\
is the identity of the group and for every }$k$\textit{\ we have }$%
S_{k}^{-1}=S_{3n-k}$.\bigskip

The proof is elementary and we omit it. We notice that for every $k$ we have 
$B\circ S_{k}=B$ and if for some meromorphic function $U$ we have $B\circ
U=B,$ then $U$ must be one of the transformations $S_{k}$ (see \cite{3}).
Consequently, the group $\{S_{k}\}$ is the group of cover transformations of 
$(\widehat{\mathbb{C}},B).$

\ 

We illustrate these results by taking $n=3$ in $(4).$
For $n=3$ and $\tau =1$ the equations $(6)$ become:\begin{equation}u\frac{r^{2}-u^{2}}{1-r^{2}u^{2}}=\omega _{k},\qquad
k=0,1,2,\end{equation}
with the obvious solution $u=-1$ for $k=0.$ We notice that for $k=1,2,$
since $\omega _{1}=\overline{\omega }_{2},$ if $u$ is a solution of one of
the two corresponding equations $(9)$ then $\overline{u}$ is a solution of
the other one, therefore these six solutions are two by two complex
conjugate. The other three solutions are $-1$ and the complex conjugate
numbers $\frac{1}{2}[1+r^{2}\pm i\sqrt{4-(1+r^{2})^{2}}]=e^{\pm i\beta },$
where $\beta =\arccos \frac{1+r^{2}}{2}$. The effect of the complex conjugation of the roots $%
u_{j}$ of the equations $(9)$, expressed in terms of the corresponding roots
of the equations $(5)$,  is the  symmetry of these roots with respect to the line $%
z(t)=e^{i\alpha }t,$ $t\in R$. The root corresponding to $u=-1$ is $%
z=-e^{i\alpha }$ and the roots corresponding to $u=e^{\pm i\beta }$ are $%
z=e^{i(\alpha \pm \beta )}.$ The dependence of these last roots on $r$
through the intermediate of $\beta $ is obvious. On the website of the project \cite{cri} we illustrate the dependence on $r$ of
all the solutions of the equation $B(z)=e^{in\alpha }.$

\bigskip

In Figure 1 we took $\ n=3,$ $r=2/3$ and $\alpha =\pi /3.$ Figure 1(a-c)
shows colored arcs connecting $\zeta _{k}$ with $a,-a,b,-b$ and $0.$ Every
color corresponds to a point $\zeta _{j}$ on the unit circle. These arcs
have been obtained by simultaneous continuation (see \cite{2}) over the real
negative half axis starting from $\zeta _{k}.$ The three pictures illustrate
the status of the continuation corresponding to three different values of $%
\tau :$ \ $\tau =1-|B(b)|,$ $\tau =1$ and $\tau =1+|B(b)|.$ The arcs converging
to $a,$ $-a$ and $0$ reach these points as $w=B(z)$ reaches $0$ (Figure
1(b)). However, the arcs converging to $b$ must have reached the respective
point\ before (Figure 1(a)), since $B(b)$ is situated between $-1$ and $0.$
When $w$ varies between $B(b)$ and $0,$ one of these arcs go from $b$ to $0$
and the other from $b$ to $a,$ while the other arcs reach, some of them $0,$
and some of them $-a.$ When $w$ reaches $0,$ all the solutions of the
equation $B(z)=-1$ are connected to some of the points $a,$ $-a,$ $b$ and $%
0 $. However, there is no arc connecting $-a$ and $0.$ We need to let $1-\tau $
vary through negative values up to $-|B(b)|$ in order to allow the pre-image
arcs starting in $-a$ and $0$ meet each other in $-b$ (Figure 1(c)). After
reflecting them into the unit circle, we obtain the fundamental domains
(Fig. 1(d)). It appears that this figure shows $0$ and $-a$ and $-1/%
\overline{a}$ as branch points of order four, since there are four arcs
converging to each one of them. But this comes in contradiction with the
expression\ $(4)$ for $n=3,$ according to which they should have order
three. However, the colors of the curves help  solve this mystery.
Indeed, in every one of these points there are curves of just three colors
meeting: red, purple and light blue in $a$ and $1/\overline{a};$ orange, blue and
magenta in $-a$ and $-1/\overline{a}$ and yellow, pink and green  in $0.$
We notice also that in $b$ and $1/\overline{b}$ \ curves of just two colors, red and green, 
are touching each other  while two curves, one blue and
the other green, are passing through $-b$ and $-1/\overline{b}.$ Therefore,
despite the appearance of $b$ and $1/\overline{b}$ being branch points of
order four and that of $-b$ and $-1/\overline{b}$ being regular points, the
reality is that all of them are branch points of order two. These facts
appear more obvious when we inspect the picture 1(d). Indeed, there the
configurations are the same in $a,$ $1/\overline{a},$ $-a,$ $-1/\overline{a}$
and $0$ as well as \ in $b,$ $-b,$ $1/\overline{b}$ and $-1/\overline{b}.$ Figure 1(g) shows the projection of the boundaries of the fundamental domains from Fig. 1(d) onto the Steiner surface. 

When extending the simultaneous continuation over the part of real half axis
from $0$ to $-|B(b)|,$ the arcs through $0,$ $a,$ and $-a$ will extend also
inside the domains already bounded by them and the unit circle. The same is
true for the domains obtained from them by reflection with respect to the
unit circle. In the end, some of the fundamental domains appear as having
inner boundaries (arcs whose points are accessible in two distinct ways from
the inside of the domain), which produce in Figure 1(f) and 1(g) the
"fractures" inside  some fundamental domains. We think that the
respective fractures are due to the fact that, in this
case,  Mathematica computes separately the values of the function on the two borders, as if indeed
they represented different parts of the boundary, and the errors of
approximation accumulate differently at the same point producing slightly
different values. There are no real inner boundaries in any fundamental
domain of this example. 

We obtain a visualization of the mapping $(4)$ by
coloring a set of annuli centered at the
origin of the $w$-plane in different colors and with saturation increasing counter-clockwise and brightness increasing outward (the saturation is determined by the argument of the point and the brightness is determined by the modulus) and imposing the same color, saturation and  brightness  to the
pre-image of every point in these annuli. 
Figures 1(e) and 1(f) show the fundamental domains of the mapping (4). The annuli are rendered in Figures 1(h) and (i). Note that in all figures we show only a selection of the annuli whose pre-images are displayed. A complete collection of annuli can be viewed on the website of the project \cite{cri}.  

\ \ \ \ \ 

If the exponent $p$ in $(3)$ is negative, or if some of $a_{k}$ have the
module greater than $1$, then $B$ is a Blaschke quotient. In this case,
parts of the interior of the unit disc switch with parts of its exterior
when mapped by $B$, yet the unit circle is still mapped on itself. However,
there might be some other curves mapped on the unit circle. This is the case
when in $(4)$ we make one of the following two changes: we replace the
factor $z^{n}$ by $1/z^{n},$ or we let $|a|$ $>1.$ First notice that
with the notation $(4)$, we have for any $a$ with $|a|$ $\neq 1:$\begin{equation} B(z)=(1/z^{n})\left[\frac{\overline{a}}{a}\frac{z^{2}-a^{2}}{\overline{a}%
^{2}z^{2}-1}\right]^{n}=B_{1/a}(1/z).\end{equation}

Therefore, when the first change is made, we need to solve an equation of
the form: $B_{1/a}(1/z)=\tau e^{-in\alpha },$ and make the same
substitution: $u=e^{-i\alpha }z.$ Instead of $(6)$ we obtain:
\begin{equation}\frac{1}{u}\frac{r^{2}-u^{2}}{1-r^{2}u^{2}}=\tau
^{1/n}\omega _{k},\qquad \tau \geq 0,\qquad k=0,1,...,n-1.\end{equation}

This time we might happen to have $B^{\prime }(\zeta )=0$ for some values $%
\zeta $ with $|\zeta |=1.$ Indeed, the equation $\ds \frac{d}{du}\left(\frac{1}{u}%
\frac{r^{2}-u^{2}}{1-r^{2}u^{2}}\right)=0$ is equivalent to:\begin{equation}r^{2}u^{4}+(1-3r^{4})u^{2}+r^{2}=0,\end{equation}
which gives:\begin{equation}u^{2}=(1/2r^{2})[3r^{4}-1\pm \sqrt{(1-r^{4})(1-9r^{4})}.\end{equation}
Then, for $1/\sqrt{3}<r<1,$ we have $\ u^{2}=e^{\pm i\gamma },$ where $%
\gamma =\arccos \frac{3r^{4}-1}{2r^{2}},$ thus $z=\pm e^{i(\alpha \pm \gamma
/2)}.$ These numbers are solutions of the equation $B^{\prime }(z)=0,$
therefore they are branch points of $(\widehat{\mathbb{C}},B)$ and they are situated
on the unit circle. This can happen only if the pre-image by\ $B$ of the
unit circle from the $w$-plane contains, besides the unit circle, two other
curves passing through these points. Due to the continuity, these must be
closed curves on $\widehat{\mathbb{C}}$ and due to the fact that $B$ commutes with $h$, they must be $h$-symmetric to each other. Moreover, the equation $%
B(z)=e^{-in\alpha }$ has all the roots $\zeta _{k},$ $k=0,1,...,3n-1$, on the
unit circle and on these two curves. Suppose that we have counted these
roots in the following way: $\zeta _{0}=-e^{i\alpha },$ then we denoted by $%
\zeta _{1}$ the root obtained after turning once counter-clockwise around
the unit circle in the $w$-plane, so that the image\ by $B$ of the arc
between $\zeta _{0}$ and $\zeta _{1}$ is that circle, etc. Finally, from $%
\zeta _{3n-1}$ we reach $\zeta _{0}$ again by turning once more around the
respective circle. In other words, the three components of the pre-image of
the unit circle from the $w$-plane can be viewed as a unique closed curve $%
\Gamma $ with self intersections in the points representing solutions
situated on the unit circle of the equation $B^{\prime }(z)=0.$ The curve $%
\Gamma $ is the lifting in the $z$-plane staring from $\zeta _{0}$ of the
closed curve obtained by tracing  the unit
circle in the $w$-plane $3n$ times counter-clockwise. When performing the simultaneous continuation over $%
w(\tau )=\tau e^{-in\alpha },\tau >0$, starting from $\zeta _{k}$, the arcs
we obtain in this way, as well as the arcs determined by $\zeta _{k}$ on $%
\Gamma $ define fundamental domains for $(\widehat{\mathbb{C}},B).$

If $r=1/\sqrt{3},$ then $\gamma =\pi $ and the self intersection points of $%
\Gamma $ are $\pm e^{i(\alpha +\pi /2)}$ on the unit circle, which are
branch points of order six for $(\widehat{\mathbb{C}},B).$ Obviously, $\Gamma $ is an 
$h$-symmetric curve (see Figure 2(b)). The sequence $(\zeta _{k})$ shows how
this curve is traced by a point $\zeta $ whose image by $B$ goes around the unit circle in the $w$-plane nine
times counter-clockwise. Figures 2(a) and 2(c) give the illusion that $\zeta $ goes
sometimes clockwise (as, for example, between $\zeta _{1}$ and $\zeta _{2},$
since we are tempted to improperly divide $\Gamma $ in connected components.  

Finally, when $r$ starts taking values less then $1/\sqrt{3},$ the curve $%
\Gamma $ separates into three disjoint components, one of which is the unit
circle and the other two are one interior to the unit circle and one
exterior to it (see Figure 2(c)). These last two components continue to
remain $h$-symmetric to each other.

In the case $r=1/\sqrt{3},$ the line determined by the points $\pm e^{i(\alpha +\pi /2)}$ intersects $%
\Gamma $ in two triplets of points having the same images by $B$ on the unit
circle in the $w$-plane. If $\zeta $ follows on $\Gamma $ and on parts of
this line the sequence $(\zeta _{k})$ as previously, we notice in all the
three cases the following. The arcs representing continuations from $\zeta _{k}$ over the real negative
half-axis meet in 14 branch points of $(\widehat{\mathbb{C}},B).$\ All these arcs
determine 18 domains which are mapped conformally by $B$ either on the open
unit disc (\emph{i-domains}), or on the exterior of the closed unit disc (%
\emph{e-domains}). Every couple of adjacent i-domains, respectively
e-domains, is separated by continuation arcs, while every couple of domains
of different types is separated by arcs of $\Gamma $. We can combine
arbitrarily two adjacent i-domains and e-domains in order to form
fundamental domains of $(\widehat{\mathbb{C}},B).$ The border of every i-domain is
obtained by traversing counter-clockwise arcs of $\Gamma $ determined by
consecutive $\zeta _{k}$ or $\zeta _{k}$ and a branch point on $\Gamma ,$ as
well as some continuation arcs between them traversed in any sense, while
for the e-domains at least one of the arcs of $\Gamma $ should be traversed
clockwise.

In the case $0<r<1/\sqrt{3}$ the equation $(13)$ has four imaginary roots $%
u_{j},$ thus the corresponding $z_{j}$ are on the line \ passing through $%
\pm e^{i(\alpha +\pi /2)}$ and they are branch points of $(\widehat{\mathbb{C}},B).$
It is easily seen that they are two by two $h$-symmetric. This line is
bordering some of i-domains and e-domains visible on Figure 2(f) and (i).
Every i-domain is mapped anticonformally by $h$ on an e-domain and
vice-versa. The dependence of the configuration on $r$ is illustrated on the website \cite{cri}. 

To summarize, Figure 2 is organized in three columns. The left column (Figure 2(a),(d), (g) and (h)) illustrates the case $1/\sqrt{3} <r<1$ and gives, respectively, the boundaries of the fundamental domains, the preimage  of colored annuli under $B$, the preimage of colored annuli in the unit disk and the projection of the boundaries of the fundamental domains on the Steiner surface. The middle column (Figure 2(b), (e), (h) and (k)) illustrates the case $r=1/\sqrt{3} $ and the right column (Figure 2(c), (f), (i) and (l)) illustrates the case  $0<r<\sqrt{3}$. Some of the colored annuli are shown in Figure 3(l).

\ 

Suppose now that\ the second change is made in $(4).$ The same conclusion as
in the previous example is valid, except that the intervals $(0,1/\sqrt{3})$
and $(1/\sqrt{3},1)$ for $\rho =1/r$ \ must be replaced by $(\sqrt{3},\infty )$%
, respectively $(1,\sqrt{3})$ for $r.$ This case is illustrated in Figure 3 and we notice a striking similarity with 
Figure 2 despite of the fact that the two cases represent
Blaschke quotients of very different nature.

\ 

If both of the previously mentioned changes are made, the formula $(10)$
applies, where $|1/a|$ $<1.$ This time $B$ switches the interior and the
exterior of the unit disc, mapping the unit circle on itself, and for $%
|\zeta |$ $=1$ we have:
\begin{equation} B^{\prime }(\zeta )=-\frac{1}{\zeta ^{2}}B_{1/a}^{\prime
}(1/\zeta )\neq 0.\end{equation}

In this case, the pre-image of the unit circle is the unit circle. Consequently,
the equation $B(z)=e^{-in\alpha }$ has distinct solutions all situated on
the unit circle. The simultaneous continuation from these points over the
ray $w(\tau )=\tau e^{-in\alpha },\tau \geq 0$, produces $3n$ arcs delimiting
the fundamental domains of $B.$ These arcs meet each other in $0,$ $a$ and $%
-a,$ which are branch points of order $n$, as well as in their symmetric
points with respect to the unit circle, and also in the four non-zero
solutions of the equation $B^{\prime }(z)=0,$ equivalent to $B_{1/a}^{\prime
}(1/z)=0.$ For $n=3$ and $a=2e^{i\pi /3},$ we obtain a Blaschke quotient. It is illustrated in Figure 4 and 
performs a similar mapping with that of Figure 1, except that the saturation of color in
Figures 4(d) and (e) is the reverse of that in Figures 1(e) and (f). It is not a surprise that the self
mapping of $P^{2}$ induced by this quotient appears to be the same as that
induced by the the Blaschke product from example 1.\bigskip \bigskip 

\section{The Case of Several Zeros of the Same Module}

Blaschke products with zeros of the same module and arguments $\alpha +2k\pi
/n$ appeared to be of special interest in \cite{1}. The condition $B\circ
h=h\circ B$ imposes the form:\begin{equation}B(z)=z^{n}\left(\frac{\overline{a}}{a}\right)^{n}\frac{z^{2n}-a^{2n}}{%
\overline{a}^{2n}z^{2n}-1},\end{equation}
where $a=re^{i\alpha }$ and $n$\ $=2k+1.$ The equation $B(z)=te^{in\alpha }$
is equivalent to $u^{3n}-r^{2n}tu^{2n}-r^{2n}u^{n}+t=0,$ where $%
u=e^{-i\alpha }z,$ which is an algebraic equation of degree $3n$ and
consequently has $3n$ solutions $u_{k}^{(j)}(t),$ $k=0,1,...,n-1,$ $j=0,1,2$
(counted with their multiplicities). In particular, $u_{k}^{(0)}(0)=0,$ $%
u_{k}^{(1)}(0)=r\omega _{k},$ and $u_{k}^{(2)}(0)=-r\omega _{k},$ $%
k=0,1,...,n-1$, while $u_{k}^{(0)}(1)=e^{-i\pi /n}\omega _{k},$ $%
u_{k}^{(1)}(1)=e^{-i\theta /n}\omega _{k}$ and $u_{k}^{(2)}(1)=e^{i\theta
/n}\omega _{k},$ $k=0,1,...,n-1$. Here $\theta $ is the argument of $%
(1+r^{2n})/2+i\sqrt{4-(1+r^{2n})^{2}}/2,$ hence $\theta =\arccos
(1+r^{2n})/2.$ Correspondingly, we have $z_{k}^{(0)}(0)=0,$ $%
z_{k}^{(1)}(0)=a\omega _{k}$ and $z_{k}^{(2)}(0)=ae^{i\pi /n}\omega _{k},$
while $z_{k}^{(0)}(1)=$ $e^{i(\alpha -\pi /n)}\omega _{k},$ $%
z_{k}^{(1)}(1)=e^{i(\alpha -\theta /n)}\omega _{k}$ and $%
z_{k}^{(2)}(1)=e^{i(\alpha +\theta /n)}\omega _{k},$ $k=0,1,...,n-1.$ Since $%
B^{\prime }(\zeta )\neq 0$ for $|\zeta |$ $=1,$ these last $3n$ solutions
are distinct.

The equation $B^{\prime }(z)=0$ has the solutions $z=0$ of order $n-1$, as
well as the simple solutions\begin{equation}b_{k}=(1/r\sqrt[2n]{2})\sqrt[2n]{3-r^{4n}-\sqrt{(3-r^{4n})^{2}-4r^{4n}}}
e^{i\alpha }\omega _{k}\end{equation}
inside the unit circle, and $1/\overline{b}_{k}$ outside the unit circle,  where $
k=0,1,...,n-1.$

As $t$ varies from $1$\ to $|B(b_{k})|,$ the points $z_{k}^{(j)}(t)$
describe $3n$ arcs expanding inside the unit circle from the points $%
z_{k}^{(j)}(1)$ situated on the unit circle. They meet in triplets in the
points $b_{k},$ when $t=|B(b_{k}))|.$ Since $|b_{k}|$ $<r,$ and $%
B(re^{i\alpha }\omega _{k})=0,$ $n$ of these arcs pass first through the
origin before arriving at $b_{k}.$ As $t$ varies from $1$ to $1/|B(b_{k})|,$
the corresponding arcs will be symmetric to the previous arcs with respect
to the unit circle. All together, they form the boundaries of $2n$ unbounded
fundamental domains and $n$ bounded fundamental domains. Let us denote by $%
\Omega _{k}^{(0)}$ the fundamental domain containing the arc $z(t)=e^{it},$ $%
t\in (\alpha +\frac{(2k-1)\pi }{n},\alpha +\frac{2k\pi -\theta }{n}),$ by $%
\Omega _{k}^{(1)}$ the fundamental domain containing the arc $z(t)=e^{it},$ $%
t\in (\alpha +\frac{2k\pi -\theta }{n},\alpha +\frac{2k\pi +\theta }{n})$ and
by $\Omega _{k}^{(2)}$ the fundamental domain containing the arc $%
z(t)=e^{it},$ $t\in (\alpha +\frac{2k\pi +\theta }{n},\alpha +\frac{%
(2k+1)\pi }{n}),$ $k=0,1,...,n-1.$ The points $0$ and $\infty $ are branch
points of order $n,$ while every point $b_{k}$ and every point $1/\overline{b%
}_{k}$ is a branch point of order $3$ of the covering Riemann surface $(%
\widehat{\mathbb{C}},B).$ \ We proved the following theorem.

\ \ 

\textbf{Theorem 3.1:} \textit{The domains }$\Omega _{k}^{(j)},$ $j=0,1,2,$ $%
k=0,1,...,n-1$ \textit{are mapped conformally by }$B$\textit{\ on }$%
\widehat{\mathbb{C}}$\textit{\ from which the the ray }$w(\tau )=\tau e^{in\alpha
},\tau \geq 0,$\textit{\ has been removed.}\bigskip

When dealing with the cover transformations of $(\widehat{\mathbb{C}},B)$ we have to
solve the equation $B(\zeta )=B(z)$  for $\zeta $, where $B$ is given by $%
(15).$ The solutions of this $3n$-degree equation are:

\begin{equation}\zeta =S_{p}^{(q)}(z)=\left[B\left|_{\Omega _{p+k \!
\!\! \!\! \pmod n}^{q+j\!\!\!\!\! \pmod
3}}\right.\right]^{-1}\circ B\left|_{\Omega _{k}^{j}}\right.(z), \ \ q,j\in \{0,1,2\}, \ \ p,k\in
\{0,1,...,n-1\}\end{equation}

It can be easily checked that $S_{p}^{(0)}(z)=\omega _{p}z,$ $%
S_{p}^{(1)}(z)=\varphi (z)\omega _{p}$ and $S_{p}^{(2)}(z)=\psi (z)\omega
_{p},$ $p=0,1,...,n-1,$ where $\varphi (z)$ and $\psi (z)$ are the uniform
branches of the following multivalued functions:\bigskip

$$z\rightarrow\left\{\left[\frac{e^{in\alpha }}{(1-(re^{-i\alpha
}z)^{2n}})\right][(r^{4n}-1)(ze^{-i\alpha })^{n}\right.\ \ \ \ \ \ \ \ \ \ \ \ \ \ \ \ \ \ \ \ \ \ \ \ \ \ \ \  \ \ \ \ \ \ \ \ \ \ \ \ \ \ $$ \begin{equation}\left.\pm \sqrt{(1-r^{4n})^{2}(ze^{-i\alpha })^{2n}+4[1-(re^{-i\alpha
}z)^{2n}][r^{2n}-(e^{-i\alpha }z)^{2n}]}]\right\}^{1/n}.\end{equation}
We notice that $\varphi (e^{i(\alpha -\pi /n)}\omega _{p})=\{e^{i(n\alpha
+\theta )}\}^{1/n}$ and $\psi (e^{i(\alpha -\pi /n)}\omega
_{p})=\{e^{i(n\alpha -\theta )}\}^{1/n},$ $p=0,1,...,n-1.$ Therefore we can
choose the principal branches of the multivalued functions $(18)$ \ such that

$$S_{p}^{(0)}(e^{i(\alpha -\pi /n)}\omega _{k})=e^{i(\alpha
-\pi /n)}\omega _{p+k},S_{p}^{(1)}(e^{i(\alpha -\pi /n)}\omega
_{k})=e^{i(\alpha +\theta /n)}\omega _{p+k}$$ 
and
$$S_{p}^{(2)}(e^{i(\alpha -\pi /n)}\omega _{k})=e^{i(\alpha
-\theta /n)}\omega _{p+k}.$$
In other words, we have: 
\begin{equation}S_{p}^{(q)}\circ S_{r}^{(s)}=S_{p+r \!\!\!\!\! \pmod n}^{q+s
\!\!\!\!\! \pmod 3}.\end{equation}
Thus, we have proved the following theorem\bigskip

\textbf{Theorem 3.2: }\emph{The group of covering transformations of }$(%
\widehat{\mathbb{C}},B),$\emph{\ with }$B$\emph{\ given by }$(15)$\emph{\ is the group
of transformations }$(17)$\emph{\ with the composition law }$(19).$\emph{\ }

\

On the website of the project \cite{cri} we illustrate the situation when $n=3.$
After the change of variable $%
\zeta =e^{-i\alpha }z$, the equation $B(z)=e^{3\alpha i}$ becomes :\begin{equation}\zeta ^{3}\frac{\zeta ^{6}-r^{6}}{r^{6}\zeta ^{6}-1}%
=1,\end{equation}
with the solutions: $\zeta _{1}=-1,$ $\zeta _{2,3}=\frac{1}{2}\pm i\frac{%
\sqrt{3}}{2},$ and the cubic roots of $\ds \frac{1}{2}[1+r^{6}\pm i\sqrt{4-(1+r^{6})^{2}}]$
for the other six solutions. We notice that the first three of them do not
depend on $r,$ while the last six tend two by two to the cubic roots of
unity as $r\rightarrow1$. Rotating the picture around the origin by an angle $\alpha $
we obtain the description of the situation in terms of\ $\ z.$ The
expression $(16)$ becomes:\begin{equation}\pm \frac{1}{r\sqrt[6]{2}}\sqrt[6]{3-r^{12}\pm \sqrt{(3-r^{12})^{2}-4r^{12}}%
}e^{i\alpha }\omega _{k},\qquad k=0,1,2,\end{equation}
where $\omega _{k}$ are the roots of order $3$ of unity. These points
are situated on the line passing through\ the origin and $a\omega _{k}$ and
they are branch points of $(\widehat{\mathbb{C}},B).$

\ 

Figure 5(a) shows the unit circle cut by 6 arcs joining $b_{k}$ and $%
1/\overline{b}_{k},$ as well as three lines passing through the origin and $%
b_{k},$ from which the segments between $b_{k}$ and $1/\overline{b}_{k}$ are
removed. There are also three rays $z(t)=te^{i[\alpha +(2k+1)\pi /3]},$ $%
k=0,1,2.$ The domains bounded by two consecutive such arcs are mapped
conformally by  $B$  on the  $w$-plane from which the positive real
half-axis has been removed.

\bigskip

In formula $(15)$ we can carry out the types of changes we performed in $(4)$ in
order to obtain Blaschke quotients commuting with $h$: replacing the factor $%
z^{n}$ by $(1/z)^{n},$ taking $|a|$ $>1,$ or both. If we denote by $%
B_{a}$ the Blaschke product $(15)$, and\begin{equation}B(z)=\left(\frac{1}{z}\right)^{n}\left(\frac{\overline{a}}{a}\right)^{n}\frac{z^{2n}-a^{2n}}{\overline{a}^{2n}z^{2n}-1},\end{equation}
then we find again that $B(z)=B_{1/a}(1/z)$ and the same arguments apply as
in the previous section. 

\bigskip

\section{Color Visualization of Blaschke Self-Mappings of the Steiner's
Surface}

The fundamental domains $\Omega _{k}$ of the Blaschke quotients studied in
the previous sections are all symmetric with respect to $h,$ i.e. $z\in
\Omega _{k}$ if and only if $h(z)\in \Omega _{k}.$ Consequently, when
factoring by the two element group $\langle h\rangle$ generated by $h,$ every
fundamental domain $\Omega _{k}$ $\subset \widehat{\mathbb{C}}$ is mapped two-to-one
on a domain $\widetilde{\Omega }_{k}\subset P^{2}.$ Moreover, the function $%
\mathbf{b}:P^{2}\rightarrow P^{2}$ defined by $\mathbf{b}(\widetilde{z})=\widetilde{%
B(z)}$ maps every $\widetilde{\Omega }_{k}$ bijectively on $P^{2},$ hence $%
\widetilde{\Omega }_{k}$ are fundamental domains of $\mathbf{b}.$

Steiner's Roman surface is a topological realization of $P^{2}=\widehat{\mathbb{C}}/ \langle h\rangle$. Endowed with the dianalytic structure induced by the analytic
structure of $\widehat{\mathbb{C}},$ the surface $P^{2}$ becomes a non orientable
Klein surface. The canonical projection
$$\Pi :\widehat{\mathbb{C}}\rightarrow P^{2}$$
defined by $\Pi (z)=\widetilde{z}$ is a morphism of Klein surfaces and $(%
\widehat{\mathbb{C}},\Pi )$ is a covering surface of $P^{2}$ (the orientable double
cover of $P^{2}).$ We have $\Pi \circ h=\Pi $ and for every Blaschke
quotient $(3),$ the identity $\Pi \circ B=\mathbf{b}\circ \Pi $ is true,
hence $\Pi \circ B\circ h=\mathbf{b}\circ \Pi \circ h=\mathbf{b}\circ \Pi .$

The mapping $\mathbf{b}$ defined on $P^{2}$ is a dianalytic mapping of the
interior of every fundamental domain $\widetilde{\Omega }_{k}$ on $P^{2}$
provided with a slit. We call it a \textit{Blaschke self-mapping} of $P^{2}.$

We have proved in \cite{4} that the Blaschke quotient $(3)$ induces a
dianalytic self-mapping of $P^{2}$ having exactly $2(p+n)+1$ fundamental
domains.

For the topological Steiner surface, the canonical image $\widetilde{%
\partial D}$ of the unit circle $\partial D$ does not play any special role.
However, when considering $P^{2}$ as a non orientable Klein surface, $%
\widetilde{\partial D}$ assumes a special role. Namely, we have shown in \cite{4} that every dianalytic self-mapping $\mathbf{b}$\ of $P^{2}$ such that $%
\mathbf{b}(\widetilde{z})\in \widetilde{\partial D}$ if and only if $%
\widetilde{z}\in \widetilde{\partial D}$ is a finite Blaschke self-mapping
of $P^{2}.$ Moreover, the boundary of every fundamental domain $\widetilde{%
\Omega }_{k}$ of $\mathbf{b}$ contains a sub-arc of $\widetilde{\partial D}.$%
\ More exactly, there are $2(n+p)+1$ distinct points $\widetilde{\zeta }%
_{k}\in \widetilde{\partial D}$ such that $\mathbf{b}(\widetilde{\zeta }%
_{k})=\widetilde{1}$ and every half open sub-arc of $\widetilde{\partial D}$
determined by two consecutive points $\widetilde{\zeta }_{k}$ and $%
\widetilde{\zeta }_{k+1}$ belongs to a unique $\widetilde{\Omega }_{k\text{ }%
}$and is mapped bijectively by $\mathbf{b}$\textbf{\ }on $\widetilde{%
\partial D},$ while $\widetilde{\Omega }_{k}$ is mapped bijectively on $%
P^{2}.$

Infinite Blaschke products $(3)$ induce infinite Blaschke self-mapping on $%
P^{2}.$ Let $E$ be the set of accumulation points of the zeros of $B.$ Since 
$B(a_{k})=0$ if and only if $B(-a_{k})=0,$ we conclude that $e^{i\theta }\in
E$ if and only if $-e^{i\theta }=h(e^{i\theta })\in E$. Therefore, there exists $%
\widetilde{E}\subset \widetilde{\partial D}$ such that $E=\Pi ^{-1}(%
\widetilde{E}).$ The function $\mathbf{b}$ cannot be defined on $\widetilde{E%
},$ because $B$ is not defined on $E.$ However, the formula $\mathbf{b}(%
\widetilde{z})=\widetilde{B(z)}$ defines $\mathbf{b}$ everywhere on $P^{2}\, \backslash \, \widetilde{E}$ and $P^{2}\, \backslash\,
\widetilde{E}=\cup _{k=1}^{\infty }\widetilde{\Omega }_{k},$ where $%
\widetilde{\Omega }_{k}$ are disjoint and $\mathbf{b}$ maps every $
\widetilde{\Omega }_{k}$ bijectively on $P^{2},$ the mapping being
dianalytic in the interior of $\widetilde{\Omega }_{k.}.$

Figures 1(g),  2(j), 2(k), 2(l), 3(j), 3(k), 4(b), and 5(b) show the domains of $P^2$ mapped bijectively on $P^2$ by  the Blaschke self-mapping of $P^{2}$ corresponding
to the Blaschke quotients defined in the corresponding previous sections. 
\bigskip

\section{Invariants of Blaschke Self-Mappings of $P^{2}$}

Let us denote by $B$ an arbitrary Blaschke quotient. If $B$ is infinite, we
denote as usual by $E$ the set of accumulation points of the zeros of $B$,
otherwise we take $E=\emptyset.$ We call an invariant of $B$ any
self-mapping $U$ of $\widehat{\mathbb{C}}\, \backslash \, E$ such that $B\circ U=B.$
Obviously, the set of invariants of $B$ is a group of transformations of $\widehat{\mathbb{C}}\, \backslash \, E$. We have proved in \cite{4} that if $B$ commutes
with $h,$ then every invariant $U$ of $B$ also commutes with $h.$ Moreover,
if $\mathbf{b}$ is the Blaschke self-mapping of $P^{2}\, \backslash\, 
\widetilde{E}$ defined by $\mathbf{b}(\widetilde{z})=\widetilde{B(z)},$ then
the\ self-mapping $u$ of $P^{2}\, \backslash \, \widetilde{E}$ defined by $u(%
\widetilde{z})=\widetilde{U(z)}$ is an invariant of $\mathbf{b},$ i.e. $%
\mathbf{b}\circ u=\mathbf{b.}$ The set of invariants of $\mathbf{b}$ is a
group of transformations of $P^{2}\, \backslash \, \widetilde{E}.$ They
represent the group of covering transformations of the Klein covering
surface $(P^{2}\, \backslash \, \widetilde{E},\mathbf{b})$ of $P^{2}.$

Once the fundamental domains $\Omega _{k}$ of $B$ are known, any invariant $%
U_{k}$ of $B$ is given by the formula $U_{k}\left|_{\Omega _{j}}\right.=[B\left|_{\Omega
_{k+j}}\right.]^{-1}\circ B\left|_{\Omega _{j}}\right.$ for every fundamental domain $\Omega _{j%
\text{ }}$of $B$. It is an easy exercise to show that composition
law is $U_{p}\circ U_{q}=U_{p+q \!\!\pmod n}$ if $\ B$ is finite of degree $%
n.$ The inverse transformations are of the form $U_{k}^{-1}=U_{n-k}$ and the
identity is $U_{0}.$ If $B$ is infinite, then $p,q\in Z$ and the composition
law is simply $U_{p}\circ U_{q}=U_{p+q}.$ The inverse transformation of $%
U_{k}$ is $U_{k}^{-1}=U_{-k}.$ The explicit computation of $U_{k}$ is in
general rather tedious. If, for example, $B$ is given by $(4),$ such a
computation would require solving equations of degree three (with literal
coefficients) of the form:$$\zeta ^{\prime }[r^{2}-\zeta ^{\prime 2}]/[1-r^{2}\zeta
^{\prime 2}]=\zeta \omega _{k}[r^{2}-\zeta ^{2}]/[1-r^{2}\zeta
^{2}],\qquad k=0,1,...,n-1,$$
where $\zeta ^{\prime }=e^{-i\alpha }z^{\prime },$ and $\zeta =e^{-i\alpha
}z.$

Suppose however that the solutions of these equations, representing the $3n$
invariants of $B$ have been found and they are of the form $z^{\prime
}=U_{k}^{(j)}(z),$ $j=0,1,2;$ $k=0,1,...,n-1.$ Then $U_{k}^{(j)}\circ
U_{k^{\prime }}^{(j^{\prime })}=U_{k^{\prime \prime }}^{(j^{\prime \prime
})},$ where $k^{\prime \prime }$ and $j^{\prime \prime }$ are uniquely
determined by $k,j,k^{\prime },j^{\prime }.$

A similar statement is true for the invariants of $B$ in $(15).$ Here, with
the notation $\zeta ^{\prime }=e^{-in\alpha }z^{\prime n}$ and $\zeta
=e^{-in\alpha }z^{n}$, we need to solve first the degree three
equations$$\zeta ^{\prime }[r^{2}-\zeta ^{\prime 2}]/[1-r^{2}\zeta
^{\prime 2}]=\zeta \lbrack r^{2}-\zeta ^{2}]/[1-r^{2}\zeta ^{2}]$$
and then replace $\zeta ^{\prime (j)}$ by $\zeta ^{(j)}\omega _{k}.$ If we
denote by $z^{\prime }=U_{k}^{(j)}(z)$ the final solutions, then it can be
checked again that the composition law above is still true.

To find the invariants of $\mathbf{b}$, we can either use the fundamental
domains $\widetilde{\Omega }_{k},$ and apply similar formulas to those used 
previously, or define directly $u_{k}$ by identities of the form $u_{k}(%
\widetilde{z})=\widetilde{U_{k}(z)}.$

\bigskip
\newpage

\begin{center} Figure 1
\end{center}

\vspace{.2in}

\noindent  \includegraphics[width=2truein,height=2truein]{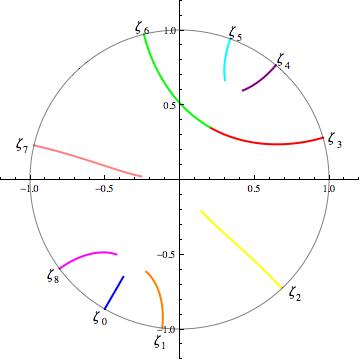} \ \ 
  \includegraphics[width=2truein,height=2truein]{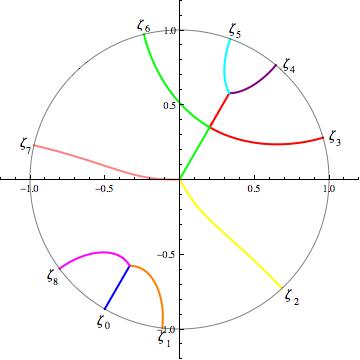} \ \  \includegraphics[width=2truein,height=2truein]{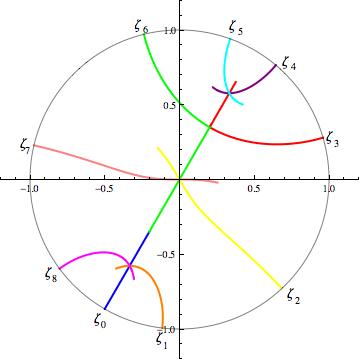}\newline
  \hspace*{.8in} $(a)$ \hspace{2in} $(b)$
   \hspace{2in} $(c)$\vspace{.1in}

  \noindent  \includegraphics[width=2truein,height=2.2truein]{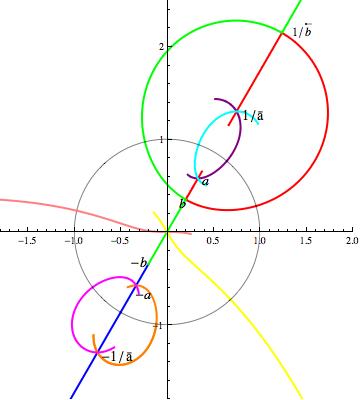} \ \  \includegraphics[width=2.3truein,height=2.3truein]{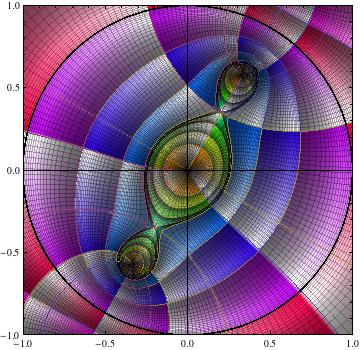} \ \ 
  \includegraphics[width=2.3truein,height=2.3truein]{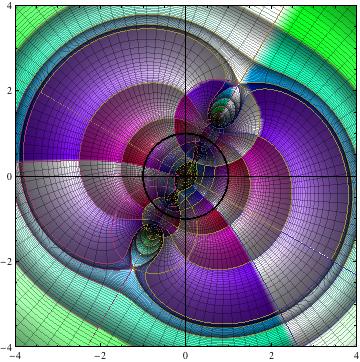}
\newline
  \hspace*{.8in} $(d)$ \hspace{2in} $(e)$
   \hspace{2in} $(f)$\vspace{.1in}
   
     \noindent   \includegraphics[width=2.5truein,height=2.8truein]{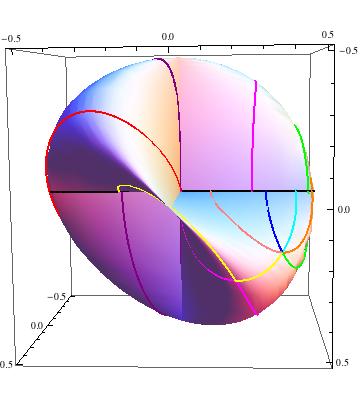} \ \  \includegraphics[width=1.95truein,height=1.8truein]{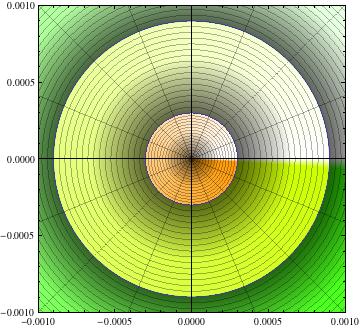}\ \ \includegraphics[width=1.8truein,height=1.8truein]{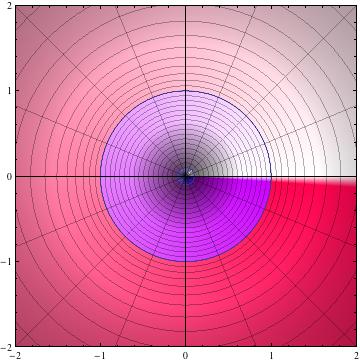}\newline
  \hspace*{1in} $(g)$ \hspace{2.2in} $(h)$  \hspace{1.7in} $(i)$\newpage

  \begin{center} Figure 2
\end{center}


\noindent  \includegraphics[width=2truein,height=2truein]{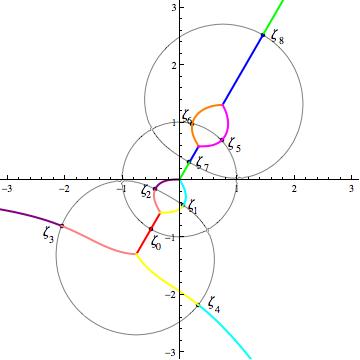} \ \ 
  \includegraphics[width=2truein,height=2truein]{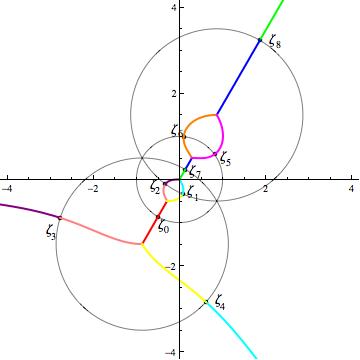} \ \  \includegraphics[width=2truein,height=2truein]{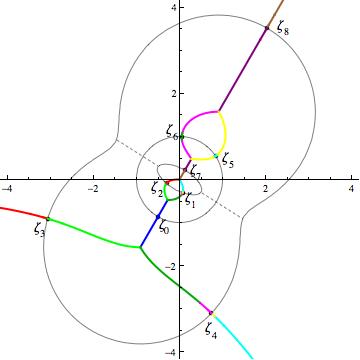}\newline
  \hspace*{.8in} $(a)$ \hspace{2in} $(b)$
   \hspace{2in} $(c)$\vspace{.05in}

     \noindent  \includegraphics[width=2truein,height=2truein]{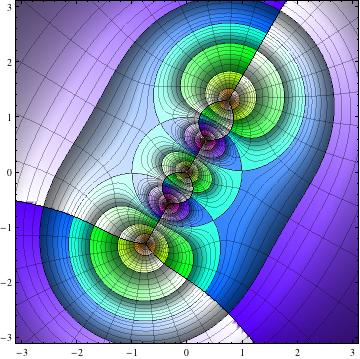} \ \  \includegraphics[width=2truein,height=2truein]{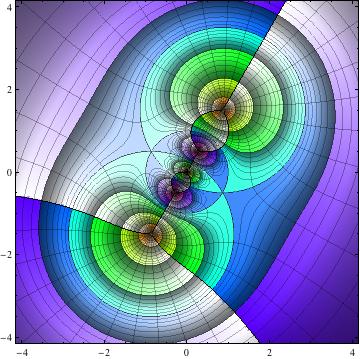} \ \ 
  \includegraphics[width=2truein,height=2truein]{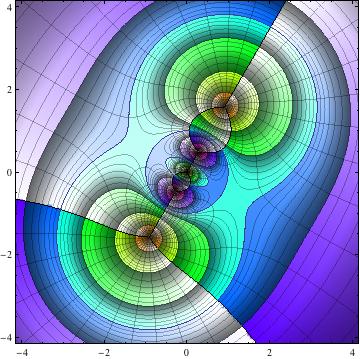}
\newline
  \hspace*{.8in} $(d)$ \hspace{2in} $(e)$
   \hspace{2in} $(f)$\vspace{.05in}

     \noindent  \includegraphics[width=2truein,height=2truein]{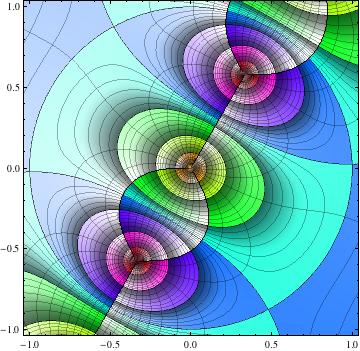} \ \  \includegraphics[width=2truein,height=2truein]{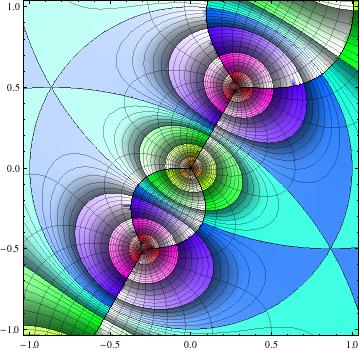} \ \ 
  \includegraphics[width=2truein,height=2truein]{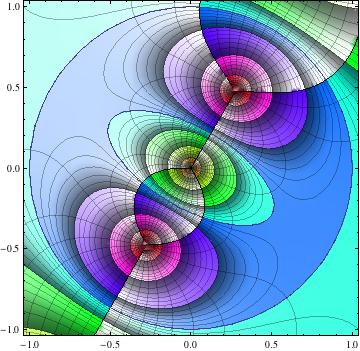}
\newline
  \hspace*{.8in} $(g)$ \hspace{2in} $(h)$
   \hspace{2in} $(i)$\vspace{.05in}

  \noindent  \includegraphics[width=2truein,height=2truein]{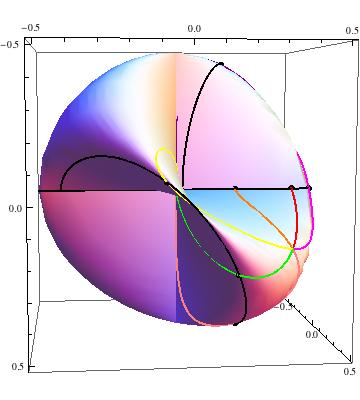} \ \  \includegraphics[width=2truein,height=2truein]{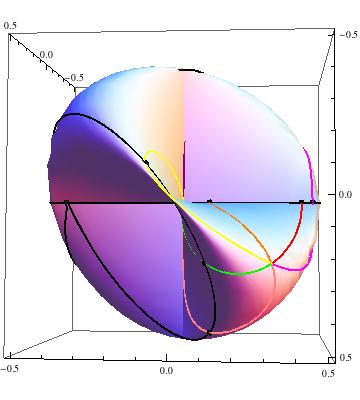} \ \ 
  \includegraphics[width=2truein,height=2truein]{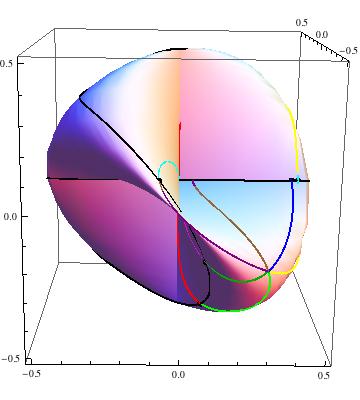}
\newline
  \hspace*{.8in} $(j)$ \hspace{2in} $(k)$
   \hspace{2in} $(l)$

\newpage

\begin{center} Figure 3
\end{center}


\noindent  \includegraphics[width=2truein,height=2truein]{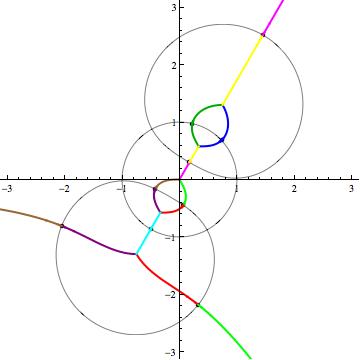} \ \ 
  \includegraphics[width=2truein,height=2truein]{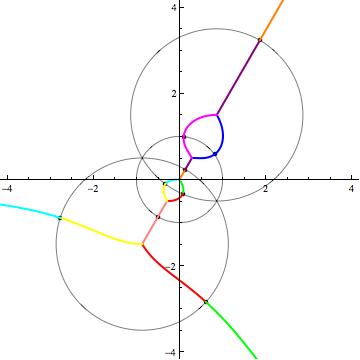} \ \  \includegraphics[width=2truein,height=2truein]{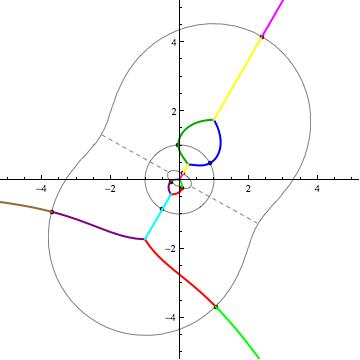}\newline
  \hspace*{.8in} $(a)$ \hspace{2in} $(b)$
   \hspace{2in} $(c)$\vspace{.05in}

     \noindent  \includegraphics[width=2truein,height=2truein]{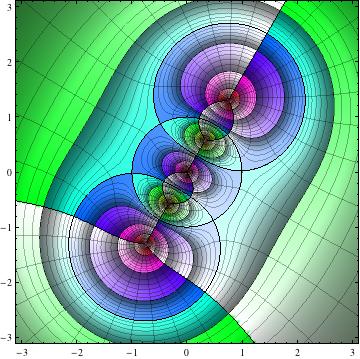} \ \  \includegraphics[width=2truein,height=2truein]{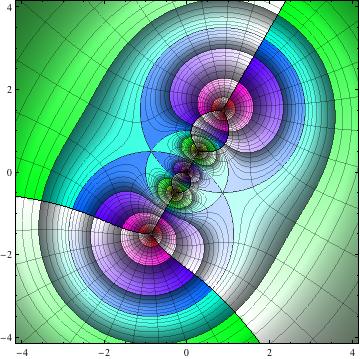} \ \ 
  \includegraphics[width=2truein,height=2truein]{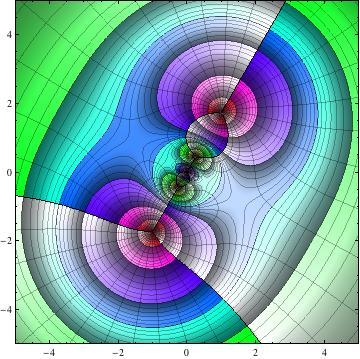}
\newline
  \hspace*{.8in} $(d)$ \hspace{2in} $(e)$
   \hspace{2in} $(f)$\vspace{.05in}

     \noindent  \includegraphics[width=2truein,height=2truein]{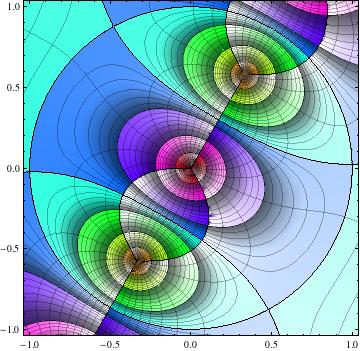} \ \  \includegraphics[width=2truein,height=2truein]{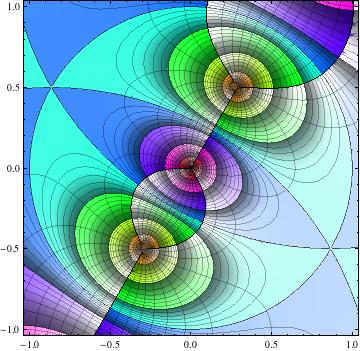} \ \ 
  \includegraphics[width=2truein,height=2truein]{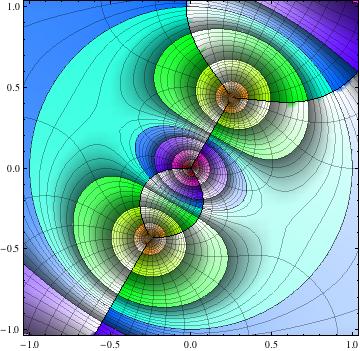}
\newline
  \hspace*{.8in} $(g)$ \hspace{2in} $(h)$
   \hspace{2in} $(i)$\vspace{.05in}

  \noindent  \includegraphics[width=2truein,height=2truein]{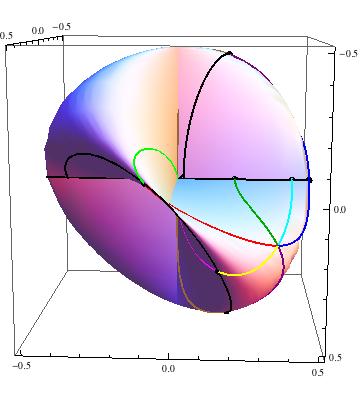} \ \  \includegraphics[width=2truein,height=2truein]{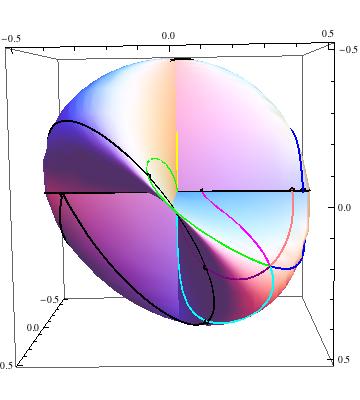} \ \ 
  \includegraphics[width=2truein,height=2truein]{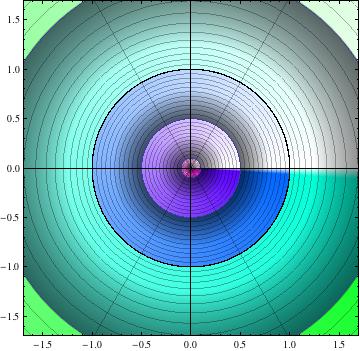}
\newline
  \hspace*{.8in} $(j)$ \hspace{2in} $(k)$
   \hspace{2in} $(l)$

\newpage

\begin{center} Figure 4
\end{center}

\includegraphics[width=2.2truein,height=2.5truein]{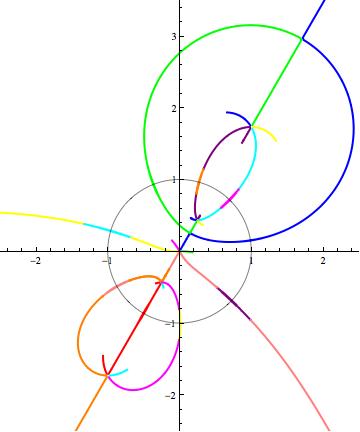}  \ \ \ \ \ \ \ \ \ \ \ \  \ \ \ \ \ \ 
\includegraphics[width=3truein,height=3truein]{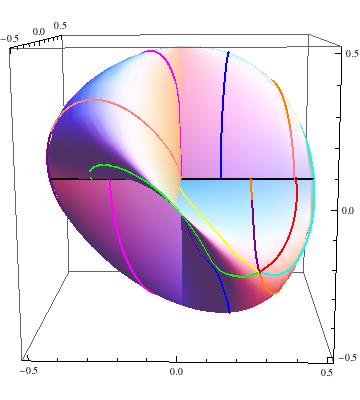}
\newline
  \hspace*{1in} $(a)$ \hspace{3.5in} $(b)$
 \vspace{.1in}

     \noindent  \includegraphics[width=3truein,height=3truein]{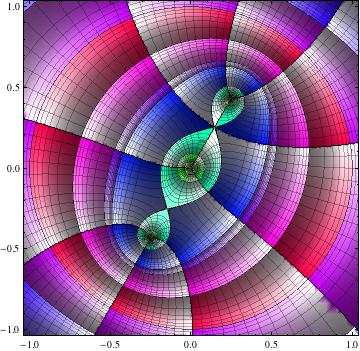} \ \ \ \ \ \ \ \ \ \ \ \ \includegraphics[width=3truein,height=3truein]{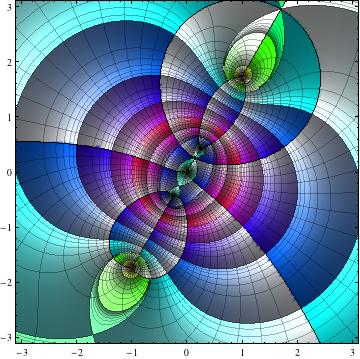} \
\newline
  \hspace*{1in} $(d)$ \hspace{3.5in} $(e)$
 \vspace{.1in}

     \noindent  \includegraphics[width=2.2truein,height=2truein]{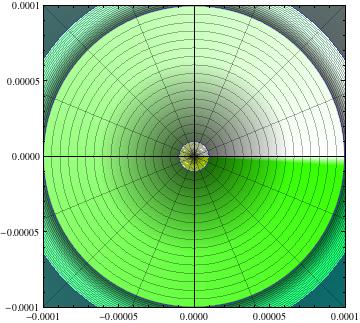}\ \ \ \ \ \ \ \ \ \ \ \  \ \ \ \ \ \ \ \ \ \ \ \  \ \ \ \ \ \ \ \ \ \ \ \includegraphics[width=2truein,height=2truein]{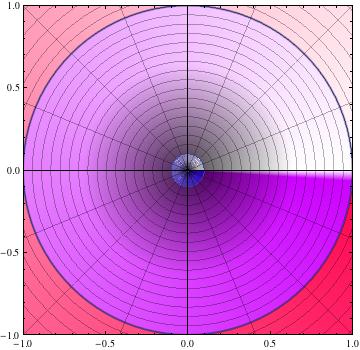} \ \ 
\newline
  \hspace*{1in} $(g)$ \hspace{3.5in} $(h)$
\vspace{.05in}

\newpage

\begin{center} Figure 5
\end{center}

\vspace{-.1in}

\noindent  \includegraphics[width=3truein,height=3truein]{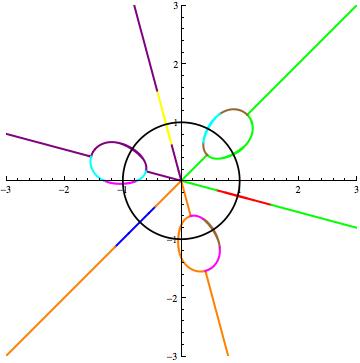} \ \ \ \ \ \ \ \ \ \ \ \
 \includegraphics[width=3truein,height=3truein]{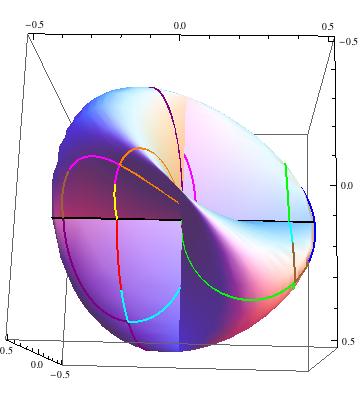}\newline
  \hspace*{1in} $(a)$ \hspace{3.5in} $(b)$
 \vspace{.2in}

     \noindent  \includegraphics[width=2.5truein,height=2.5truein]{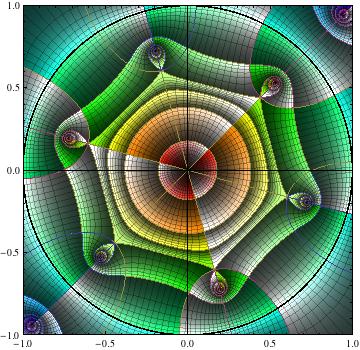} \ \ \ \ \ \ \ \ \ \ \ \ \ \ \ \ \ \ \ \ \includegraphics[width=2.5truein,height=2.5truein]{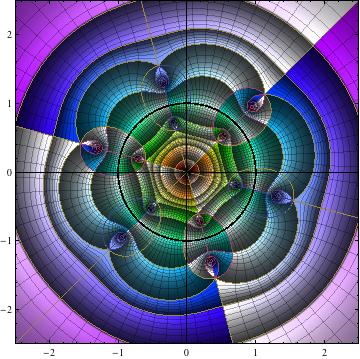} \
\newline
  \hspace*{1in} $(d)$ \hspace{3.5in} $(e)$
 \vspace{.2in}

     \noindent  \includegraphics[width=2.5truein,height=2.5truein]{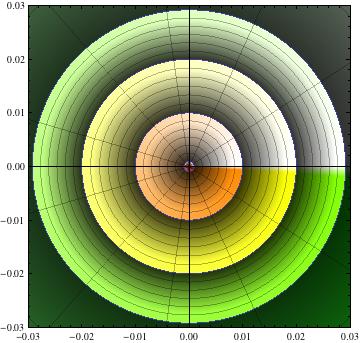}\ \ \ \ \ \ \ \ \ \ \ \ \ \ \ \ \ \ \ \ \ \includegraphics[width=2.5truein,height=2.5truein]{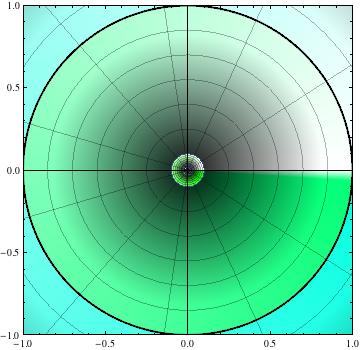} \ \ 
\newline
  \hspace*{1in} $(g)$ \hspace{3.5in} $(h)$
\vspace{.05in}
\newpage

\section{Technical details} 

All images for this article have been created using the software  \emph{Mathematica 6} on a MacBookPro with a 2.33 GHz Intel Core 2 Duo processor. The formulas used for projecting curves on the Steiner surface follow \cite{welke}. Sample code is available on the website of the project \cite{cri}. 

\section*{Acknowledgements}

The authors would like to thank Szabolcs Horv\'at for  helpful hints on Mathematica.

\bigskip 

\bigskip 

\noindent College of the Holy Cross, Worcester, Massachusetts, USA, cballant@holycross.edu

\ 

\noindent York University, Toronto, Ontario, Canada, dghisa@yorku.ca

\end{document}